\documentclass[final,10pt]{asme2ej}

\usepackage{amsmath}         
\usepackage{amssymb}
\usepackage{mathrsfs}
\usepackage{graphicx}
\usepackage{color}
\usepackage{enumerate}
\usepackage{multirow}
\usepackage{multicol}
\usepackage{booktabs}
\usepackage{epstopdf}

\newtheorem{remark}{Remark}                             \newtheorem{algorithm}{Algorithm}
\newtheorem{example}{Example}

\title{Consistent approximation of fractional order operators}

\author{Yiheng Wei
    \affiliation{
	School of Mathematics, \\
    Southeast University,\\
    Nanjing 211189, China\\
    Email: neudawei@ustc.edu.cn\\
    }	
}
\author{YangQuan Chen
    \affiliation{
	School of Engineering,\\
    University of California, Merced,\\
    CA 95343, USA\\
    Email: ychen53@ucmerced.edu
    }	
}
\author{Yingdong Wei
    \affiliation{
	Department of Automation, \\
    University of Science and Technology of China,\\
    Hefei, 230026, China\\
    Email: kb8@mail.ustc.edu.cn
    }	
}
\author{Xuefeng Zhang
    \affiliation{
	School of Sciences,\\
    Northeastern University, \\
    Shenyang 110819, China\\
    Email: zhangxuefeng@mail.neu.edu.cn
    }	
}

\begin{document}

\maketitle

\begin{abstract}
{\it
Fractional order controllers become increasingly popular due to their versatility and superiority in various performance. However, the bottleneck in deploying these tools in practice is related to their analog or numerical implementation. Numerical approximations are usually employed in which the approximation of fractional differintegrator is the foundation. Generally, the following three identical equations always hold, i.e., $\frac{1}{s^\alpha}\frac{1}{s^{1-\alpha}} = \frac{1}{s}$, $s^\alpha \frac{1}{s^\alpha} = 1$ and $s^\alpha s^{1-\alpha} = s$. However, for the approximate models of fractional differintegrator $s^\alpha$, $\alpha\in(-1,0)\cup(0,1)$, there usually exist some conflicts on the mentioned equations, which might enlarge the approximation error or even cause fallacies in multiple orders occasion. To overcome the conflicts, this brief develops a piecewise approximate model and provides two procedures for designing the model parameters. The comparison with several existing methods shows that the proposed methods do not only satisfy the equalities but also achieve high approximation accuracy. From this, it is believed that this work can serve for simulation and realization of fractional order controllers more friendly.
}
\end{abstract}


\section{Introduction}\label{Section 1}
Fractional calculus, as a generalization of the classical calculus of integrals and derivatives, has gained considerable popularity during the past three decades, mainly due to its demonstrated applications in numerous seemingly diverse and widespread fields of science and engineering \cite{Sun:2018CNSNS}. Indeed, it does provide several potentially powerful tools for practical applications, such as fractional memristor \cite{Chen:2019NN}, fractional capacitor \cite{Semary:2019JAR}, fractional inductor \cite{Adhikary:2018TCSI}, fractional resonator \cite{Adhikary:2016TCSI}, fractional oscillator \cite{Alejandro:2020JAR}, fractional controllers \cite{Zhang:2017ISA,Zhang:2020ISA,Chen:2020TCSII}, fractional models \cite{Shi:2020VRa,Shi:2020VRb}. As for the recent progress on the related topics, the readers can refer to the mentioned excellent papers and the references therein.

The numerical approximation of fractional differintegrator $s^\alpha$ becomes crucial with the increasing demand of applications on fractional order systems. Many methods were developed while none has a clear overwhelming advantage \cite{Vinagre:2000FCAA}. This situation always existed until a seminal recursive algorithm was proposed by Oustaloup \cite{Oustaloup:2000TCSI}. This algorithm provided a procedure for hardware implementation not merely established a numerical approximation method. Now, it is widely used in practice. Inspired by this heuristic algorithm, a similar algorithm was proposed for the fractional integrator instead of the fractional differentiator \cite{Poinot:2003SP}. Afterwards, a series of effective and efficient algorithms were established \cite{Krishna:2011SP,Meng:2012DSMC,Romero:2013ISA,Wei:2014IJCAS,Pakhira:2015ISA,Wei:2016ISA,Tavazoei:2016JAS}. For example, a model with summation form was derived for the fractional differentiator, which is convenient for circuit realization \cite{Abdelaty:2018TCSII}. After deriving the approximate model for the fractional differintegrator, the stability issue was discussed in \cite{Deniz:2016ISA,Sabatier:2018Alg}. A direct discretization method producing low integer order discrete-time transfer functions was proposed for fractional order transfer functions \cite{De:2018ISA}. A block diagram based simulation scheme was developed for fractional order systems under Caputo definition \cite{Bai:2018ISA}. To make the degree of the approximate model low without the loss of approximation accuracy, the fixed-pole approximation scheme was established \cite{Liang:2014IJSS,Liang:2017JCSC,Wei:2018ISA}. To enhance the robustness of the approximate model, the multiple pole method was introduced \cite{Li:2020ISA,Wei:2021DSMC}.


Despite the great achievements, they have a major drawback which is the invalidation of the related identical equations. In active design, this implies bringing accumulative error with multiple fractional orders and greatly reduces the practicability. For example, when we consider fractional internal model control, fractional non-overshoot control and fractional oscillation circuit design, the related identical equations become more and more important. For the convenience of narration, we set ${\mathscr I}^\alpha(s)=\frac{1}{s^\alpha}$ and ${\mathscr D}^\alpha(s)={s^\alpha}$. The corresponding approximate operators are ${\hat{ \mathscr I}_\kappa^\alpha }\left( s \right)$ and ${\hat{ \mathscr D}_\kappa^\alpha }\left( s \right)$, respectively, where $\kappa$ is the index to distinguish the different approximations. Motivated by the above discussion, the objective of this work is to design the alternative, accurate, available and attractive approximation algorithms such that the following conditions are satisfied.
\begin{enumerate}[i)]
    \item ${\hat{ \mathscr I}_\kappa^\alpha }\left( s \right){\hat{ \mathscr I}_\kappa^{1 - \alpha }}\left( s \right) = \frac{1}{s}$;
    \item ${\hat{ \mathscr D}_\kappa^\alpha }\left( s \right){\hat { \mathscr I}_\kappa^{\alpha }}\left( s \right) = 1$;
    \item ${\hat{ \mathscr D}_\kappa^\alpha }\left( s \right){\hat{ \mathscr D}_\kappa^{1 - \alpha }}\left( s \right) = s$.
\end{enumerate}

Contributions of this work are mainly manifested in three aspects. (a) The conflict on the order of  fractional approximation models is discussed for the first time. To be more exact, they are wrong in the mentioned occasion. (b) A piecewise approximation model is proposed to avoid the conflict, which extends the application range and improves the approximation accuracy. (c) Two independent procedures for designing the model parameters are developed and discussed, which greatly facilitate the design and implementation of fractional differintegrators. (d) The multiple pole principle is adopted to enrich the applicability.

Bearing this in mind, Section \ref{Section 2} proposes a universal framework to approximate ${\mathscr I}^\alpha(s)$ and then four sets of parameters are developed, resulting in ${\hat {\mathscr I}_\kappa^\alpha }\left( s \right)$ and ${\hat {\mathscr D}_\kappa^\alpha }\left( s \right)$, $\kappa=1,2,3,4$. Section \ref{Section 3} provides two illustrative examples to show the effectiveness of the proposed methods, demonstrating its added value with respect to the state of art. Finally, some concluding remarks are drawn in Section \ref{Section 4}.

\section{Main Results}\label{Section 2}
\subsection{Basic framework}
In the beginning, let us approximate the fractional integrator ${\mathscr I}^\alpha(s)$ with the piecewise model
\begin{equation}\label{Eq1}
{\textstyle {\hat {\mathscr I}_\kappa^\alpha }\left( s \right) = \left\{ \begin{array}{l}
K\prod\nolimits_{i = 1}^n {{{\big( {\frac{{s + {z_i}}}{{s + {p _i}}}} \big)}^k}} ,0 < \alpha  \le 0.5,\\
\frac{{\bar K}}{s}\prod\nolimits_{i = 1}^n {{{\big( {\frac{{s + {{\bar z}_i}}}{{s + {{\bar p}_i}}}} \big)}^k},0.5 < \alpha  < 1,}
\end{array} \right.}
\end{equation}
where $\kappa=1,2,3,4$, $k\in\mathbb{Z}_+$, $n\in\mathbb{Z}_+$, $K$, $\bar K$ are the gains, $-p_i$, $-\bar p_i$ are the poles and $-z_i$, $-\bar z_i$ are the zeros. For this model, when the following relationships
\begin{equation}\label{Eq2}
{\textstyle
\left\{
\begin{array}{rl}
{\bar p_i} =&\hspace{-0pt} \left.{z _i} \right|_{\alpha = 1 - \alpha},\\
{\bar z_i} =&\hspace{-0pt} \left. {p_i} \right|_{\alpha = 1 - \alpha},\\
\bar K =&\hspace{-0pt} {{ K^{-1} }} |_{\alpha = 1 - \alpha},
\end{array}\right.}
\end{equation}
are satisfied, the condition i) holds for $\alpha\in(0,0.5)\cup(0.5,1)$.

To meet the condition ii) for $\alpha\in(0,1)$, the approximate model for fractional differentiator ${\mathscr D}^\alpha(s)$ is designed as
\begin{equation}\label{Eq3}
{\textstyle {\hat {\mathscr D}_\kappa^\alpha }\left( s \right) =\frac{1}{{\hat {\mathscr I}_\kappa^\alpha }\left( s \right)}.}
\end{equation}
With this design, the condition iii) is also implied.

Up to now, all the elaborated conditions have been met. We will focus on how to design the parameters subsequently. Since ${\hat {\mathscr D}_\kappa^\alpha }\left( s \right)$ can be uniquely identified with the given ${\hat {\mathscr I}_\kappa^\alpha }\left( s \right)$, only the design of ${\hat {\mathscr I}_\kappa^\alpha }\left( s \right)$ will be discussed hereafter, i.e., $z_i$, $p_i$, $i=1,2,\cdots,n$ and $K$.


Given the parameters $\kappa, k, n$ and the interested frequency interval $[\omega_l,\omega_h]$, it is expected that $p_i,z_i,\bar p_i,\bar z_i\in[\omega_l,\omega_h]$ for any $i=1,2,\cdots,n$. First, for the case of $0<\alpha\le0.5$, supposing that $\hat {\mathscr I}_k^\alpha(s)$ and ${\mathscr I}^\alpha(s)$ have the same gain at the middle frequency $s={\rm j}\omega_m$, i.e.,
\begin{equation}\label{Eq4}
{\textstyle | {\hat {\mathscr I}_\kappa^\alpha ({\rm{j}}{\omega _m})} | =  | {{{\mathscr I}^\alpha }({\rm{j}}{\omega _m})} |, }
\end{equation}
the desired gain $K = | {\rm{j}}{\omega _m}|^{-\alpha}\prod\nolimits_{i = 1}^n \big|{{{ {\frac{{{\rm{j}}{\omega _m} + {p _i}}}{{{\rm{j}}{\omega _m} + {z_i}}}} }}}  \big|^k$ can be determined, where ${\omega _m} \triangleq \sqrt {{\omega _l}{\omega _h}} $, ${\rm j} \triangleq \sqrt {-1} $. By using formula (\ref{Eq2}), one has
\begin{equation}\label{Eq5}
{\textstyle \begin{array}{rl}
\bar K =&\hspace{-0pt} {K^{ - 1}}{|_{\alpha  = 1 - \alpha }}\\
 =&\hspace{-0pt} { {\big|{{({\rm{j}}{\omega _m})}^\alpha }\prod\nolimits_{i = 1}^n {{{\big(\frac{{{\rm{j}}{\omega _m} + {z_i}}}{{{\rm{j}}{\omega _m} + {p _i}}}\big)}^k}} \big|} _{\alpha  = 1 - \alpha }}\\
 =&\hspace{-0pt} \big|{({\rm{j}}{\omega _m})^{1 - \alpha }}\prod\nolimits_{i = 1}^n {{{\big(\frac{{{\rm{j}}{\omega _m} + {z_i}{|_{\alpha  = 1 - \alpha }}}}{{{\rm{j}}{\omega _m} + {p _i}{|_{\alpha  = 1 - \alpha }}}}\big)}^k}} \big|\\
 =&\hspace{-0pt} \big|{({\rm{j}}{\omega _m})^{1 - \alpha }}\prod\nolimits_{i = 1}^n {{{(\frac{{{\rm{j}}{\omega _m} + {\bar p_i}}}{{{\rm{j}}{\omega _m} + {\bar z_i}}})}^k}} \big|\\
 =&\hspace{-0pt}|{\rm{j}}{\omega _m}|^{1 - \alpha }\prod\nolimits_{i = 1}^n \big|\frac{{{\rm{j}}{\omega _m} + {\bar p_i}}}{{{\rm{j}}{\omega _m} + {{\bar z}_i}}}\big|^k,
\end{array} }
\end{equation}
which means that (\ref{Eq4})also holds for the $0.5<\alpha<1$ case.

Next, two independent design procedures will be provided for $p_i$, $z_i$, $\bar p_i$, $\bar z_i$, $i=1,2,\cdots,n$.

\subsection{Two-point boundary value based method}
\begin{figure}[!htbp]
  \centering
  \includegraphics[width=0.85\textwidth]{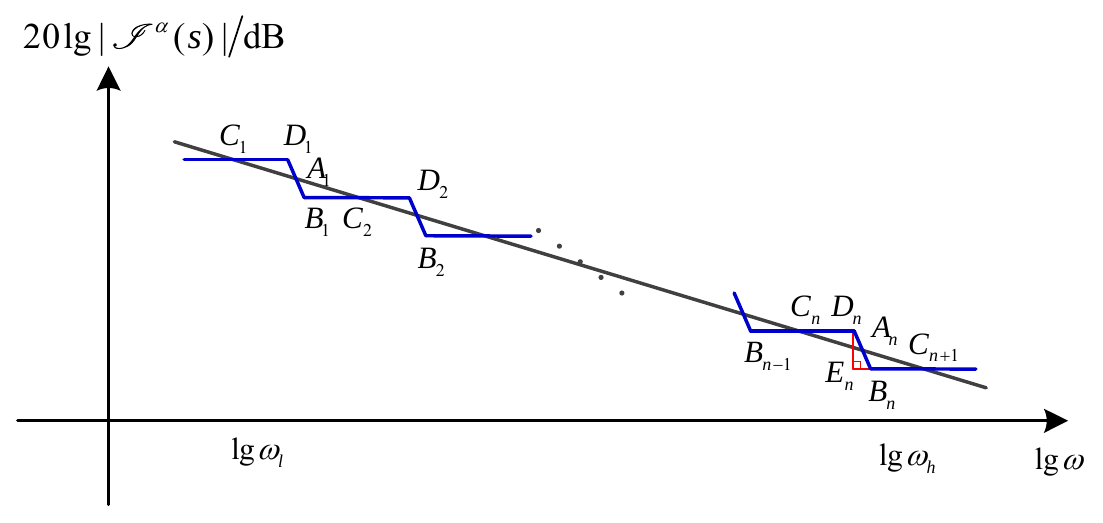}
  \caption{Integrator approximation in the Bode diagram. (The $0<\alpha\le0.5$ case)}\label{Fig 1}
\end{figure}

On the principle of the poles and the zeros alternately appear $k$ times by $k$ times, the asymptotic lines of the Bode magnitude diagram for ${{ \mathscr I}^\alpha }\left( s \right)$ and ${\hat{ \mathscr I}_\kappa^\alpha }\left( s \right)$ with $0<\alpha\le0.5$ are shown in Fig. \ref{Fig 1}.

For Fig. \ref{Fig 1}, it is noted that $B_i$ and $D_i$ correspond to the zeros and the poles and their abscissas are ${\rm lg}\,z_i$ and ${\rm lg}\,p_i$, respectively. For convenience, let us define the abscissa of points $A_i$, $C_i$ as ${\rm lg}\,\rho_i$, ${\rm lg}\,\sigma_i$, $i=1,2,\cdots,n$. The slope of line $C_1C_{n+1}$ is $-20\alpha\,{\rm dB/dec}$. The slope of line $C_iD_i$ is $ 0\,{\rm dB/dec}$ and the slope of line $D_iB_i$ is $ -20 k\,{\rm dB/dec}$. Since $0<\alpha<k$, along this way, the approximation problem can be solved via approximating a line with the slope $-20\alpha\,{\rm dB/dec}$ by a set of lines with the slopes $ 0\,{\rm dB/dec}$ and $ -20 k\,{\rm dB/dec}$ within the logarithmic frequency range $[{\rm lg}\,\omega_l,{\rm lg}\,\omega_h]$. The approximate error can be described by the total area of all triangle formed by composite curves and the target line.

As shown in the previous study \cite{Wei:2018ISA}, when the triangles are congruent, that is,
\begin{equation}\label{Eq6}
{\textstyle \bigtriangleup{A_i}{D_i}{C_i}\cong\bigtriangleup{A_{i + 1}}{D_{i + 1}}{C_{i + 1}}\cong\bigtriangleup{A_{i}}{B_i}{C_{i+ 1}}},
\end{equation}
$i=1,2,\cdots,n$, the approximate error is minimum in the sense of the area. From this relationship, one has
\begin{equation}\label{Eq7}
{\textstyle {\rm lg}\,\rho_i-{\rm lg}\,p_{i}={\rm lg}\,z_{i}-{\rm lg}\,\rho_i},
\end{equation}
\begin{equation}\label{Eq8}
{\textstyle 2{\rm lg}\,\rho_i={\rm lg}\,\sigma_{i}+{\rm lg}\,\sigma_{i+1}}.
\end{equation}

Assume that Point $E_n$ is the intersection point of Line $D_nE_n$ and Line $B_n C_{n+1}$ such that Line $D_nE_n$ is vertical. Since both $B_n E_n$ and $C_{n+1}E_n$ are horizontal lines, the slope information implies
\begin{equation}\label{Eq9}
{\textstyle \left\{ \begin{array}{rl}
\frac{{|{D_n}{E_n}|}}{{|{B_n}{E_n}|}} =&\hspace{-0pt} 20k,\\
\frac{{|{D_n}{E_n}|}}{{|{C_{n + 1}}{E_n}|}} =&\hspace{-0pt} 20\alpha ,
\end{array} \right.}
\end{equation}
which shows
\begin{equation}\label{Eq10}
{\textstyle {\rm lg}\,\rho_{i}-{\rm lg}\,p_i=\frac{\alpha}{k}({\rm lg}\,\rho_{i}-{\rm lg}\, \sigma_{i})}.
\end{equation}

When the crossing points $C_1$, $C_{n+1}$ are the boundary points, i.e., ${\sigma _1} = {\omega _l}$, ${\sigma _{n+1}} = {\omega _h}$, then ${\sigma _i}={\omega _l}{\big( {\frac{{{\omega _h}}}{{{\omega _l}}}} \big)^{\frac{{i - 1}}{n}}}$. By combining formulas (\ref{Eq7}), (\ref{Eq8}) and (\ref{Eq10}), one is ready to obtain the value of ${p _i}$, ${z _i}$. According to (\ref{Eq2}), $\bar p_i,\bar z_i$ can be calculated immediately. Additionally, it is proven that the desired condition $p_i,z_i,\bar p_i,\bar z_i\in[\omega_l,\omega_h]$ always holds for any $k\in\mathbb{Z}_+$, $\alpha\in(0,1)$ and $i=1,2,\cdots,n$. From this, a novel approximate model can be summarized as Algorithm \ref{Algorithm 1}.

\begin{algorithm}\label{Algorithm 1}
The approximate model ${\hat {\mathscr I}_1^\alpha }\left( s \right)$ for the fractional integrator ${\mathscr I}^\alpha(s)$ with $0<\alpha<1$ can be designed in the form of (\ref{Eq1}) with ${p _i}= {\omega _l}{\big( {\frac{{{\omega _h}}}{{{\omega _l}}}} \big)^{\frac{{2i - 1 - \alpha /k}}{{2n}}}}$, ${z _i} = {\omega _l}{\big( {\frac{{{\omega _h}}}{{{\omega _l}}}} \big)^{\frac{{2i - 1 + \alpha /k}}{{2n}}}}$, ${\bar p_i} = {\omega _l}{\big( {\frac{{{\omega _h}}}{{{\omega _l}}}} \big)^{\frac{{2i - 1 + 1/k - \alpha /k}}{{2n}}}}$, $
{{\bar z}_i}={\omega _l}{\big( {\frac{{{\omega _h}}}{{{\omega _l}}}} \big)^{\frac{{2i - 1 - 1/k + \alpha /k}}{{2n}}}}$, $i=1,2,\cdots,n$, $K = | {\rm{j}}{\omega _m}|^{-\alpha} \prod\nolimits_{i = 1}^n \big|  {\frac{{{\rm{j}}{\omega _m} + {p _i}}}{{{\rm{j}}{\omega _m} + {{z }_i}}}} \big|^k$ and $\bar K = |{\rm{j}}{\omega _m}|^{1 - \alpha }\prod\nolimits_{i = 1}^n \big|\frac{{{\rm{j}}{\omega _m} + {\bar p_i}}}{{{\rm{j}}{\omega _m} + {{\bar z}_i}}}\big|^k$.
\end{algorithm}

When the turning points $D_1$, $B_n$ are the boundary points, i.e., ${p _1} = {\omega _l}$, ${z _{n}} = {\omega _h}$, combining formulas (\ref{Eq7}), (\ref{Eq8}) with (\ref{Eq10}), the desired value of ${p _i}$, ${z _i}$ can be obtained. In this case, a new approximate model is expressed as Algorithm \ref{Algorithm 2}.

\begin{algorithm}\label{Algorithm 2}
The approximate model ${\hat {\mathscr I}_2^\alpha }\left( s \right)$ for the fractional integrator ${\mathscr I}^\alpha(s)$ with $0<\alpha<1$ can be designed in the form of (\ref{Eq1}) with ${p _i} = {\omega _l}{\big( {\frac{{{\omega _h}}}{{{\omega _l}}}} \big)^{\frac{{i - 1}}{{n - 1{\rm{ + }}\alpha /k}}}}$, ${{z_i}}= {\omega _l}{\big( {\frac{{{\omega _h}}}{{{\omega _l}}}} \big)^{\frac{{i - 1 + \alpha /k}}{{n - 1{\rm{ + }}\alpha /k}}}}$, ${{\bar p_i}}= {\omega _l}{\big( {\frac{{{\omega _h}}}{{{\omega _l}}}} \big)^{\frac{{i - 1 + 1/k - \alpha /k}}{{n - 1 + 1/k - \alpha /k}}}}$, ${\bar z_i} = {\omega _l}{\big( {\frac{{{\omega _h}}}{{{\omega _l}}}} \big)^{\frac{{i - 1}}{{n - 1 + 1/k - \alpha /k}}}}$, $i=1,2,\cdots,n$, $K = | {\rm{j}}{\omega _m}|^{-\alpha} \prod\nolimits_{i = 1}^n \big|  {\frac{{{\rm{j}}{\omega _m} + {p _i}}}{{{\rm{j}}{\omega _m} + {{z }_i}}}} \big|^k$ and $\bar K = |{\rm{j}}{\omega _m}|^{1 - \alpha }\prod\nolimits_{i = 1}^n \big|\frac{{{\rm{j}}{\omega _m} + {\bar p_i}}}{{{\rm{j}}{\omega _m} + {{\bar z}_i}}}\big|^k$.
\end{algorithm}

In Algorithms \ref{Algorithm 1} and \ref{Algorithm 2}, the parameters $p_i,z_i$ are designed firstly. Then, the parameters $\bar p_i,\bar z_i$ can be calculated via formula (\ref{Eq2}). In subsequent discussion, we will consider the $0.5<\alpha<1$ case firstly, and then derive the model for the $0<\alpha\le0.5$ case.

\begin{figure}[!htbp]
  \centering
  \includegraphics[width=0.85\textwidth]{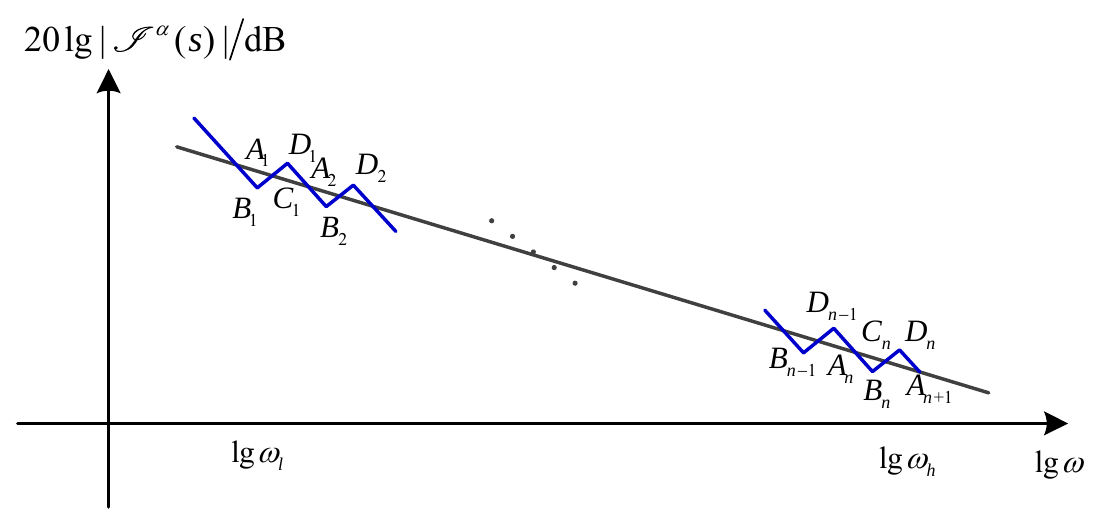}
  \caption{Integrator approximation in the Bode diagram. (The $0.5<\alpha<1$ case)}\label{Fig 2}
\end{figure}

With the given conditions, the asymptotic lines of the Bode magnitude diagram for ${{ \mathscr I}^\alpha }\left( s \right)$ and ${\hat{ \mathscr I}_\kappa^\alpha }\left( s \right)$ with $0.5<\alpha<1$ are shown in Fig. \ref{Fig 2}.

To achieve better approximation performance, the following condition is assumed.
\begin{equation}\label{Eq11}
{\textstyle \bigtriangleup{A_i}{B_i}{C_i}\cong\bigtriangleup{A_{i + 1}}{B_{i + 1}}{C_{i + 1}}\cong\bigtriangleup{A_{i+ 1}}{D_i}{C_{i}}}.
\end{equation}

Similarly, define the abscissa of points $A_i,B_i,C_i$ and $D_i$ as ${\rm lg}\,\varrho_i$, ${\rm lg}\,\bar z_i$, ${\rm lg}\,\delta_i$, and ${\rm lg}\bar p_i$, for any $i=1,2,\cdots,n$. The slope of line $C_1C_{n+1}$ is $-20\alpha\,{\rm dB/dec}$. The slope of line $A_iB_i$ is $ -20\,{\rm dB/dec}$ and the slope of line $B_iD_i$ is $ 20(k-1)\,{\rm dB/dec}$. Then one has the following relations
\begin{equation}\label{Eq12}
{\textstyle {\rm lg}\,\delta_i-{\rm lg}\,\bar z_{i}={\rm lg}\bar p_{i}-{\rm lg}\,\delta_i},
\end{equation}
\begin{equation}\label{Eq13}
{\textstyle 2{\rm lg}\,\delta_i={\rm lg}\,\varrho_{i}+{\rm lg}\,\varrho_{i+1}},
\end{equation}
\begin{eqnarray}\label{Eq14}
{\textstyle
\begin{array}{rl}
2\left( {\lg {{\bar z }_i} - \lg {{\varrho }_i}} \right)=&\alpha \left( {\lg {{\varrho }_{i + 1}} - \lg {{\varrho }_i}} \right)\\
&+ \left( {k - 1} \right)\left( {\lg \bar {p}_i - \lg {{\bar z }_i}} \right).
\end{array}}
\end{eqnarray}

Assuming the crossing points $A_1$, $A_{n+1}$ are the boundary points, i.e., ${\varrho _1} = {\omega _l}$, ${\varrho _{n+1}} = {\omega _h}$, one has ${\varrho _i}={\omega _l}{\big( {\frac{{{\omega _h}}}{{{\omega _l}}}} \big)^{\frac{{i - 1}}{n}}}$. From (\ref{Eq12})-(\ref{Eq14}), it follows ${\bar p_i} = {\omega _l}{\big( {\frac{{{\omega _h}}}{{{\omega _l}}}} \big)^{\frac{{2i -1+1/k- \alpha /k}}{{2n}}}}$ and ${\bar z_i} = {\omega _l}{\big( {\frac{{{\omega _h}}}{{{\omega _l}}}} \big)^{\frac{{2i - 1-1/k + \alpha /k}}{{2n}}}}$. With the obtained $\bar p_i$ and $\bar z_i$, the parameters $p_i$ and $z_i$ can be derived via the rule in (\ref{Eq2}). From the continuity of the model on the order, the parameters for the $\alpha=0.5$ case can also be developed. Till now, one approximation model same to Algorithm \ref{Algorithm 1} follows. When we assume the turning points $B_1$, $D_n$ are the boundary points, i.e., ${\bar z _1} = {\omega _l}$, ${ \bar p _{n}} = {\omega _h}$, using formulas (\ref{Eq12})-(\ref{Eq14}) yields ${{\bar p_i}}= {\omega _l}{\big( {\frac{{{\omega _h}}}{{{\omega _l}}}} \big)^{\frac{{i - 1 + 1/k - \alpha /k}}{{n - 1 + 1/k - \alpha /k}}}}$ and ${\bar z_i} = {\omega _l}{\big( {\frac{{{\omega _h}}}{{{\omega _l}}}} \big)^{\frac{{i - 1}}{{n - 1 + 1/k - \alpha /k}}}}$ which are the same with Algorithm \ref{Algorithm 2}. Interestingly, different means reach the same approximation model.

\subsection{One-point limited distance based method}
\begin{figure}[!htbp]
  \centering
  \includegraphics[width=0.85\textwidth]{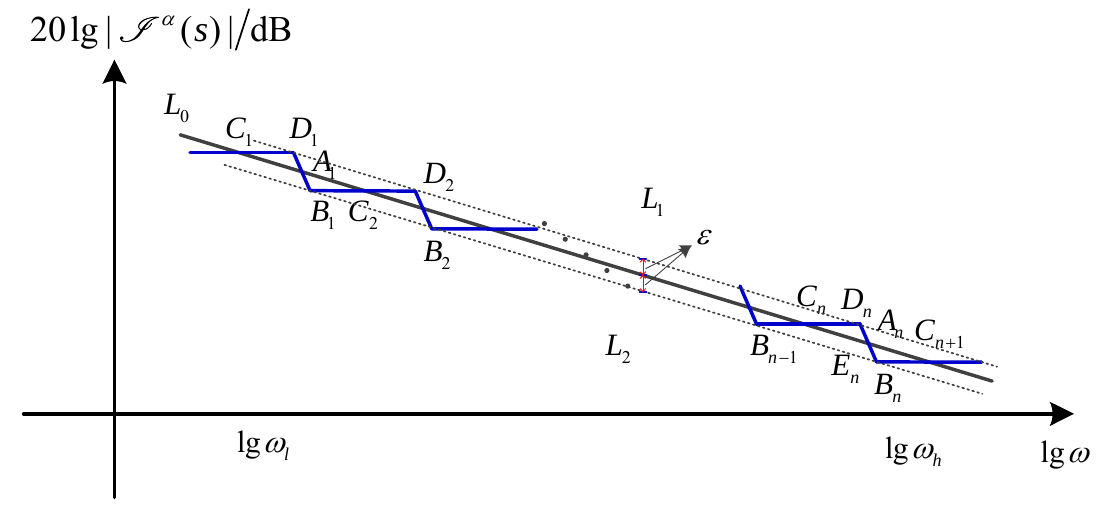}
  \caption{Integrator approximation in the Bode diagram. (The $0<\alpha\le0.5$ case)}\label{Fig 3}
\end{figure}

Hereafter, we introduce three auxiliary lines $L_0$, $L_1$ and $L_2$, which are respectively the magnitude frequency characteristic curves of $\frac{1}{s^\alpha}$, $\frac{k_1}{s^\alpha}$ and $\frac{k_2}{s^\alpha}$. We still start from the $0<\alpha\le0.5$ case. The asymptotic lines of the Bode magnitude are shown in Fig. \ref{Fig 3}, where Points $A_1$, $A_2$, $\cdots$, $A_n$ and $C_1$, $C_2$, $\cdots$, $C_{n+1}$ are on Line $L_0$, Points $D_1$, $D_2$, $\cdots$, $D_n$ are on Line $L_1$ and Points $B_1$, $B_2$, $\cdots$, $B_n$ are on Line $L_2$.

Suppose that the vertical distance between $L_0$ and $L_1$/$L_2$ is $\varepsilon>0$. Define three functions as
\begin{equation}\label{Eq15}
{\textstyle
{L_0}(\omega ) =  - 20\alpha \lg \omega ,}
\end{equation}
\begin{equation}\label{Eq16}
{\textstyle
{L_1}(\omega ) = 20\lg {k_1} - 20\alpha \lg \omega  ,}
\end{equation}
\begin{equation}\label{Eq17}
{\textstyle
{L_2}(\omega ) = 20\lg {k_2} - 20\alpha \lg \omega ,}
\end{equation}
where $\omega\in(0,+\infty)$. Based on the condition ${L_1}(\omega )-{L_0}(\omega )={L_0}(\omega )-{L_2}(\omega )=\varepsilon$, one obtains $k_1=10^{\frac{\varepsilon}{20}}$ and $k_2=10^{-\frac{\varepsilon}{20}}$.

Define the magnitude frequency characteristic curves of $\frac{h_i}{s^k}$ as $M_i(\omega )=20\lg {h_i} - 20 k\lg \omega$, $i=1,2,\cdots,n$. Assume lines $M_i$ pass through Points $D_i$ and $B_i$. From the geometrical relationship, one has
\begin{equation}\label{Eq18}
{\textstyle
\left\{\begin{array}{rl}
{L_0}({\sigma _i}) =&\hspace{-0pt} {L_1}({p _i})={L_2}({z _{i-1}}),\\
{M_i}({p _i}) =&\hspace{-0pt} {L_1}({{p }_i}),\\
{M_i}({z _i}) =&\hspace{-0pt} {L_2}({z_i}).
\end{array}\right.}
\end{equation}

Assume that the left boundary point is the crossing point $C_1$ and the right boundary point locates between $B_nB_{n+1}$, i.e., $\sigma_1=\omega_l$, $z_n\le\omega_h<z_{n+1}$. By applying formulas (\ref{Eq2}) and (\ref{Eq18}), the desired parameters are calculated. Then, the approximate model can be described as follows.

\begin{algorithm}\label{Algorithm 3}
The approximate model ${\hat {\mathscr I}_3^\alpha }\left( s \right)$ for the fractional integrator ${\mathscr I}^\alpha(s)$ with $0<\alpha<1$ can be designed in the form of (\ref{Eq1}) with ${p _i} = {10^{\frac{{\varepsilon (2ki - k - \alpha )}}{{20\alpha (k - \alpha )}}}}{\omega _l}$, ${z_i} = {10^{\frac{{\varepsilon (2ki - k + \alpha )}}{{20\alpha (k - \alpha )}}}}{\omega _l}$, ${\bar p_i} = {10^{\frac{{\varepsilon (2ki - k + 1 - \alpha )}}{{20(1 - \alpha )(k - 1 + \alpha )}}}}{\omega _l}$, ${{\bar z}_i} = {10^{\frac{{\varepsilon (2ki - k - 1 + \alpha )}}{{20(1 - \alpha )(k - 1 + \alpha )}}}}{\omega _l}$, $i=1,2,\cdots,n$, $\frac{{20\nu (k - \nu )}}{{2kn + k + \nu }}\lg \big( {\frac{{{\omega _h}}}{{{\omega _l}}}} \big) < \varepsilon  \le \frac{{20\nu (k - \nu )}}{{2kn - k + \nu }}\lg \big( {\frac{{{\omega _h}}}{{{\omega _l}}}} \big)$, $\nu = 0.5 - |\alpha  - 0.5|$, 
$K = | {\rm{j}}{\omega _m}|^{-\alpha} \prod\nolimits_{i = 1}^n \big|  {\frac{{{\rm{j}}{\omega _m} + {p _i}}}{{{\rm{j}}{\omega _m} + {{z }_i}}}} \big|^k$ and $\bar K = |{\rm{j}}{\omega _m}|^{1 - \alpha }\prod\nolimits_{i = 1}^n \big|\frac{{{\rm{j}}{\omega _m} + {\bar p_i}}}{{{\rm{j}}{\omega _m} + {{\bar z}_i}}}\big|^k$.
\end{algorithm}

Assume that the left boundary point is the turning point $D_1$ and the right boundary point locates between $B_nB_{n+1}$, i.e., $p_1=\omega_l$, $z_n\le\omega_h<z_{n+1}$. Similarly, the aforementioned conditions could give a new approximate model.

\begin{algorithm}\label{Algorithm 4}
The approximate model ${\hat {\mathscr I}_4^\alpha }\left( s \right)$ for the fractional integrator ${\mathscr I}^\alpha(s)$ with $0<\alpha<1$ can be designed in the form of (\ref{Eq1}) with ${p _i} = {10^{\frac{{\varepsilon (ki - k)}}{{10\alpha (k - \alpha )}}}}{\omega _l}$, ${z_i} = {10^{\frac{{\varepsilon (ki - k + \alpha )}}{{10\alpha (k - \alpha )}}}}{\omega _l}$, ${\bar p_i} = {10^{\frac{{\varepsilon (ki - k + 1 - \alpha )}}{{10(1 - \alpha )(k - 1 + \alpha )}}}}{\omega _l}$, ${{\bar z}_i} = {10^{\frac{{\varepsilon (ki - k)}}{{10(1 - \alpha )(k - 1 + \alpha )}}}}{\omega _l}$, $i=1,2,\cdots,n$, $\frac{{10\nu (k - \nu )}}{{kn + \nu }}\lg \big( {\frac{{{\omega _h}}}{{{\omega _l}}}} \big) < \varepsilon  \le \frac{{10\nu (k - \nu )}}{{kn - k + \nu }}\lg \big( {\frac{{{\omega _h}}}{{{\omega _l}}}} \big)$, $\nu  = 0.5 - (\alpha  - 0.5){\mathop{\rm sgn}} (\alpha  - 0.5)$, $K = | {\rm{j}}{\omega _m}|^{-\alpha} \prod\nolimits_{i = 1}^n \big|  {\frac{{{\rm{j}}{\omega _m} + {p _i}}}{{{\rm{j}}{\omega _m} + {z_i}}}} \big|^k$ and $\bar K = |{\rm{j}}{\omega _m}|^{1 - \alpha }\prod\nolimits_{i = 1}^n \big|\frac{{{\rm{j}}{\omega _m} + {\bar p_i}}}{{{\rm{j}}{\omega _m} + {{\bar z}_i}}}\big|^k$.
\end{algorithm}

\begin{remark}\label{Remark 1}
In Algorithms \ref{Algorithm 1} - \ref{Algorithm 4}, four piecewise models ${\hat {\mathscr I}_\kappa^\alpha }\left( s \right)$, $\kappa=1,2,3,4$ are constructed to approximate the fractional integrator ${\mathscr I}^\alpha(s)$ with $0<\alpha<1$. It is ready to check that as $n\to +\infty$, $\omega_l\to0$, $\omega_h\to+\infty$ and $\varepsilon\to0$, ${ {\mathscr I}^\alpha }\left( s \right)$ can be approximated to any degree of accuracy. 
By adopting the principle in (\ref{Eq3}), ${\hat {\mathscr D}_\kappa^\alpha }\left( s \right)$, $\kappa=1,2,3,4$ can be developed accordingly. To make a fair and full comparison with the existing literature, six approximate models are borrowed, i.e., $\hat {\mathscr I}_5^\alpha \left( s \right)$ from \cite{Poinot:2003SP}, $\hat {\mathscr D}_5^\alpha \left( s \right)$ from \cite{Oustaloup:2000TCSI}, $\hat {\mathscr I}_6^\alpha \left( s \right)$ and $\hat {\mathscr D}_6^\alpha \left( s \right)$ from \cite{Li:2020ISA}, $\hat {\mathscr I}_7^\alpha \left( s \right)$ and $\hat {\mathscr D}_7^\alpha \left( s \right)$ from \cite{Sabatier:2018Alg}, respectively. For better replicate and verify the provided results, the related models are listed as follows.
\begin{equation}\label{Eq19}
{\textstyle {\hat {\mathscr I}_5^\alpha }\left( s \right) = \frac{K}{s}\prod\nolimits_{i = 1}^n {\frac{{s + {{z}_i}}}{{s + {p _i}}}} ,}
\end{equation}
\begin{equation}\label{Eq20}
{\textstyle {\hat {\mathscr D}_5^\alpha }\left( s \right)= \bar K\prod\nolimits_{i = 1}^n {\frac{{s + {{\bar z}_i}}}{{s + {\bar p_i}}}} ,}
\end{equation}
where ${p _i} = {\omega _l}{\big( {\frac{{{\omega _h}}}{{{\omega _l}}}} \big)^{\frac{{i - \alpha }}{{n - \alpha }}}}$, ${{z }_i} = {\omega _l}{\big( {\frac{{{\omega _h}}}{{{\omega _l}}}} \big)^{\frac{{i - 1}}{{n - \alpha }}}}$, ${\bar p_i} = {\omega _l}{\big( {\frac{{{\omega _h}}}{{{\omega _l}}}} \big)^{\frac{{2i - 1 + \alpha }}{{2n}}}}$, ${{\bar z}_i} = {\omega _l}{\big( {\frac{{{\omega _h}}}{{{\omega _l}}}} \big)^{\frac{{2i - 1 - \alpha }}{{2n}}}}$, $i = 1,2, \cdots n$, $K = \prod\nolimits_{i = 1}^n {\big| {\frac{{{\rm{j}} + {p_i}}}{{{\rm{j}} + {{z}_i}}}} \big|} $ and $\bar K = \omega _h^\alpha$.
\begin{equation}\label{Eq21}
{\textstyle {\hat {\mathscr I}_6^\alpha }\left( s \right) = \frac{K}{s}{\prod\nolimits_{i = 1}^n {\big( {\frac{{s + {{z }_i}}}{{s + {p _i}}}} \big)} ^2} ,}
\end{equation}
\begin{equation}\label{Eq22}
{\textstyle {\hat {\mathscr D}_6^\alpha }\left( s \right)= \bar K\prod\nolimits_{i = 1}^n {{{\big( {\frac{{s + {{\bar z}_i}}}{{s + {\bar p_i}}}} \big)}^2}} ,}
\end{equation}
where ${p _i} = {\omega _l}{\big( {\frac{{{\omega _h}}}{{{\omega _l}}}} \big)^{\frac{{4i - 1 - \alpha }}{{4n}}}}$, ${{z }_i} = {\omega _l}{\big( {\frac{{{\omega _h}}}{{{\omega _l}}}} \big)^{\frac{{4i - 3 + \alpha }}{{4n}}}}$, ${\bar p_i} = {\omega _l}{\big( {\frac{{{\omega _h}}}{{{\omega _l}}}} \big)^{\frac{{4i - 2 + \alpha }}{{4n}}}}$, ${{\bar z}_i} = {\omega _l}{\big( {\frac{{{\omega _h}}}{{{\omega _l}}}} \big)^{\frac{{4i - 2 - \alpha }}{{4n}}}}$, $i = 1,2, \cdots n$, $K = \omega _h^{1 - \alpha }$, and $\bar K = \omega _h^\alpha$.
\begin{equation}\label{Eq23}
{\textstyle {\hat {\mathscr I}_7^\alpha }\left( s \right) = K\prod\nolimits_{i = 1}^n {\frac{{s + {{z }_i}}}{{s + {p _i}}}} ,}
\end{equation}
\begin{equation}\label{Eq24}
{\textstyle {\hat {\mathscr D}_7^\alpha }\left( s \right)= \frac{1}{K}\prod\nolimits_{i = 1}^n {\frac{{s + {{p }_i}}}{{s + {z _i}}}} ,}
\end{equation}
where ${p _i} = {\omega _l}{\big( {\frac{{{\omega _h}}}{{{\omega _l}}}} \big)^{\frac{{2i - 1 - \alpha }}{{2n}}}}$, ${{z }_i} = {\omega _l}{\big( {\frac{{{\omega _h}}}{{{\omega _l}}}} \big)^{\frac{{2i - 1 + \alpha }}{{2n}}}}$, $i = 1,2, \cdots n$ and $K = \prod\nolimits_{i = 1}^n {\big| {\frac{{{\rm{j}} + {{z }_i}}}{{{\rm{j}} + {p _i}}}} \big|} $.

Then the expected relations are checked in the following 7 cases and the results are shown as Table \ref{Table 1}. It is safe to say that the proposed algorithms have a distinct superiority from this point of view. The salient feature is just what we are pursuing in this work.
\begin{table}[!htbp]
\caption{The associativity for approximation operators with $\alpha\in(0,0.5)\cup(0.5,1)$.}\label{Table 1}
\centering
\scriptsize
\begin{tabular}{lccc}
\toprule
{}&condition i)&condition ii)&condition iii)\\
\hline
case 1: $\kappa=1$&$\checkmark$&$\checkmark$&$\checkmark$\\
case 2: $\kappa=2$&$\checkmark$&$\checkmark$&$\checkmark$\\
case 3: $\kappa=3$&$\checkmark$&$\checkmark$&$\checkmark$\\
case 4: $\kappa=4$&$\checkmark$&$\checkmark$&$\checkmark$\\
case 5: $\kappa=5$&$\times$&$\times$&$\times$\\
case 6: $\kappa=6$&$\times$&$\times$&$\times$\\
case 7: $\kappa=7$&$\times$&$\checkmark$&$\times$\\
\bottomrule
\end{tabular}
\end{table}
\end{remark}

\begin{remark}\label{Remark 2}
In general, with large $n$, the smaller $\epsilon$, the higher approximation accuracy. Setting $\omega_l=10^{-3}$, $\omega_h=10^3$, $n=10$ and $k=2$, Fig. \ref{Fig 4} shows the available $\epsilon$ for Algorithms \ref{Algorithm 3} and \ref{Algorithm 4}. When $\epsilon  = \frac{{10\nu (k - \nu )}}{{kn}}\lg \big( {\frac{{{\omega _h}}}{{{\omega _l}}}} \big)$, one has $\sigma_{n+1}=\omega_h$ and in this case Algorithm \ref{Algorithm 3} reduces to Algorithm \ref{Algorithm 1}. When $\epsilon  = \frac{{10\nu (k - \nu )}}{{kn - k + \nu }}\lg \big( {\frac{{{\omega _h}}}{{{\omega _l}}}} \big)$, one has $\bar \omega_{n}=\omega_h$. At this point, Algorithm \ref{Algorithm 4} degenerates into Algorithm \ref{Algorithm 2}. The degenerating condition for $\epsilon$ is marked as the special value in Fig. \ref{Fig 4}. Because there exists a suitable range for $\epsilon$ instead of a fixed value, Algorithms \ref{Algorithm 3} and \ref{Algorithm 4} could enjoy much design freedom.
\begin{figure}[!htbp]
  \centering
  \includegraphics[width=0.5\textwidth]{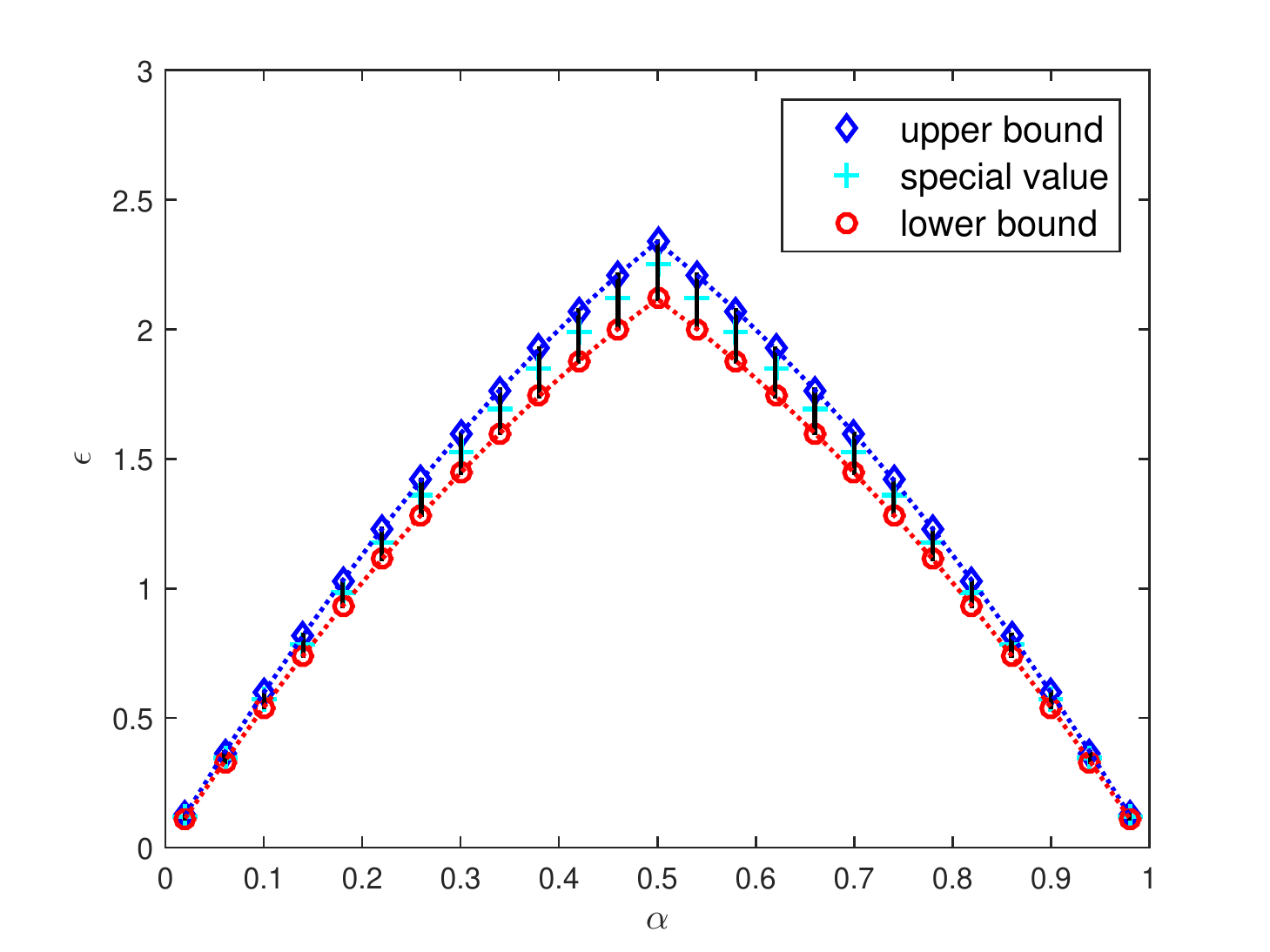}\includegraphics[width=0.5\textwidth]{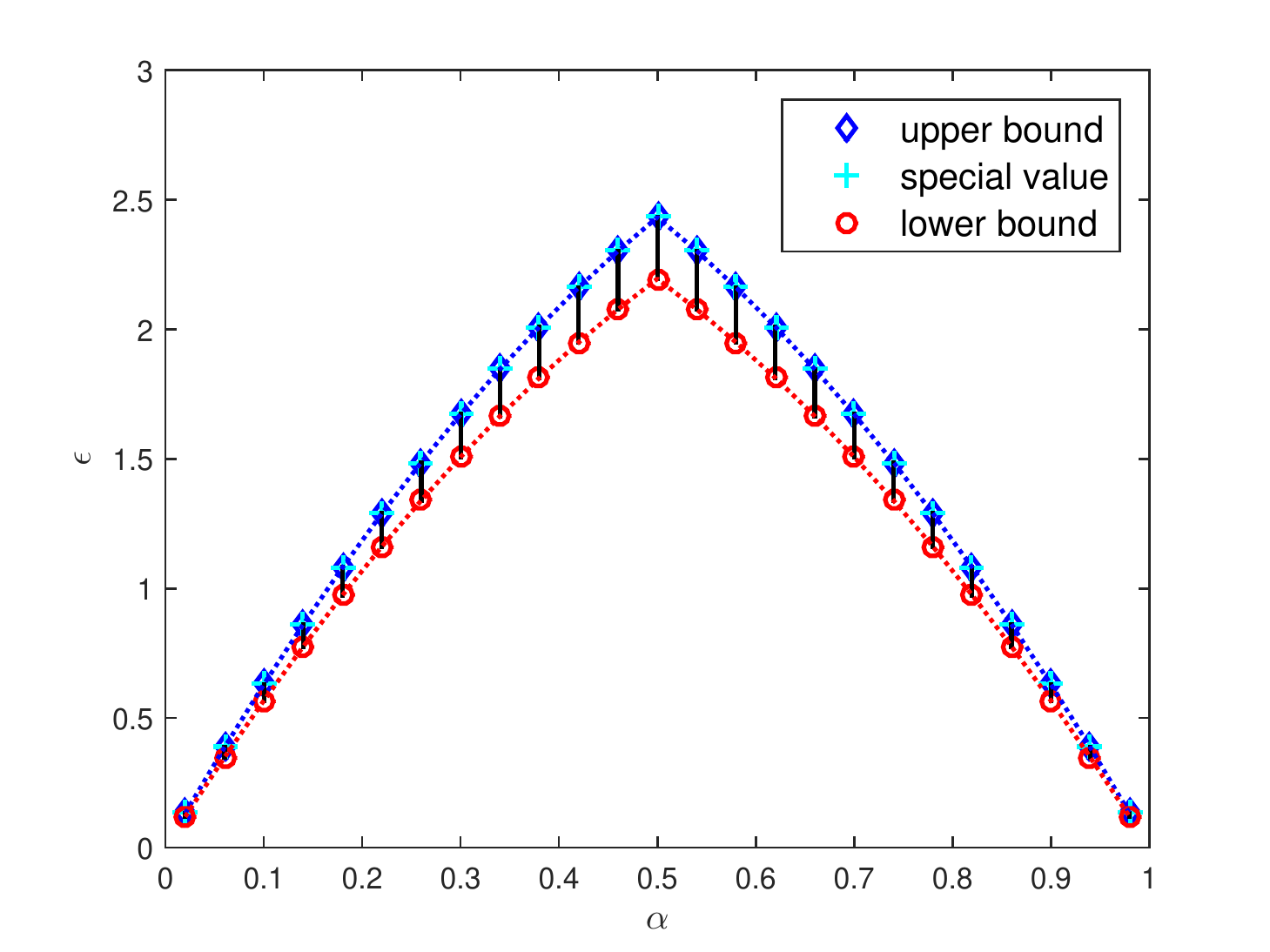}
  \caption{The range of value for $\epsilon$ with different $\alpha$. (Left: $\epsilon$ in  ${\hat {\mathscr I}_3^\alpha }( s )$; Right: $\epsilon$ in  ${\hat {\mathscr I}_4^\alpha }( s )$)}\label{Fig 4}
\end{figure}
\end{remark}

\begin{remark}\label{Remark 3}
The developed approximate model in (\ref{Eq1}) is the multiplicative case. For the single pole case, i.e., $k=1$, (\ref{Eq1}) is equivalently expressed as
\begin{equation}\label{Eq25}
{\textstyle {\hat {\mathscr I}_\kappa^\alpha }\left( s \right) = \left\{ \begin{array}{ll}
K+\sum\nolimits_{i = 1}^n {\frac{{{r_i}}}{{s + {p _i}}}}  &\hspace{-0pt},0 < \alpha  \le 0.5,\\
\frac{{{c_0}}}{s} + \sum\nolimits_{i = 1}^n {\frac{{{c_i}}}{{s + {\bar p_i}}}} &\hspace{-0pt},0.5 < \alpha  < 1,
\end{array} \right.}
\end{equation}
where ${r_i} = K\prod\nolimits_{l = 1,l \ne i}^n {\frac{{{z_l} - {p _i}}}{{{p _l} - {p _i}}}} $, ${c_0} = \bar K\prod\nolimits_{i = 1}^n {\frac{{{{\bar z}_i}}}{{{\bar p_i}}}}$, ${c_i} =  - \frac{{\bar K}}{{{\bar p_i}}}\prod\nolimits_{l = 1,l \ne i}^n {\frac{{{{\bar z}_l} - {\bar p_i}}}{{{\bar p_l} - {\bar p_i}}}} $. With this summation form, we can construct the RC circuits implementation in Fig. \ref{Fig 5} with specially selected $R_i$ and $C_i$.
\begin{figure}[!htbp]
  \centering
  \includegraphics[width=0.85\textwidth]{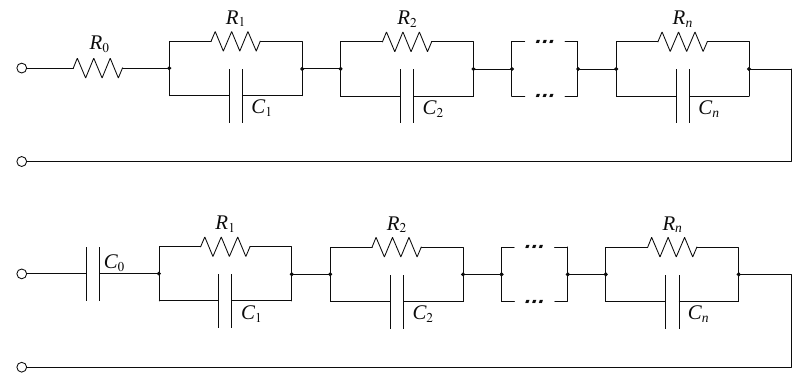}
  \caption{The RC circuits implementation for approximation fractional integrator ${\hat {\mathscr I}_\kappa^\alpha }\left( s \right)$, $\kappa=1,2,3,4$. (Above: $0<\alpha\le0.5$; Below: $0.5<\alpha<1$)}\label{Fig 5}
\end{figure}
For the multiple pole case, i.e., $k\neq1$, the equivalent summation form is derived as
\begin{eqnarray}\label{Eq26}
{\textstyle {\hat {\mathscr I}_\kappa^\alpha }\left( s \right) = \left\{ \begin{array}{ll}
K+\sum\nolimits_{i = 1}^n\sum\nolimits_{l = 1}^k {\frac{{{r_{il}}}}{{(s + {p _i})^l}}}  &\hspace{-0pt},0 < \alpha  \le 0.5,\\
\frac{{{c_0}}}{s} + \sum\nolimits_{i = 1}^n\sum\nolimits_{l = 1}^k {\frac{{{c_{il}}}}{{(s + {\bar p_i})^l}}} &\hspace{-0pt},0.5 < \alpha  < 1,
\end{array} \right.}
\end{eqnarray}
which is a little more complex than (\ref{Eq25}), while it can also be implemented by electric circuit. To sum up, in spite of the difficulties, the multiplicative approximation can be transferred into the summation form and its circuit implementation issue is not a burden.
\end{remark}

\section{Numerical Simulations}\label{Section 3}
In this section, two examples are provided to evaluate the advantages and the applicability of the developed methods.  For comparison, the mentioned seven cases in Remark \ref{Remark 1} will be considered once again. The interested frequency range $[\omega_l,\omega_h]$ $=[10^{-3},10^{3}]$, the number of zeros $n=10$, the sampling period $h=0.001$s and the number of times $k=2$ are preset. The distance variable $\varepsilon$ in Algorithms 3 and 4 are designed as the special value in Remark \ref{Remark 2}, respectively. All the related code can be found in {\tt www.mathworks.com} with title ``Fractional Differintegrator Approximation''.

\begin{example}Performance in frequency domain. 

Given the interested order $\alpha=0.1,0.2,\cdots,0.9$ and the frequency points $\omega(k) =\big(\frac{\omega_h}{\omega_l}\big)^{\frac{k-1}{10000-1}}$, $k=1,2,\cdots,10000$, the Bode diagrams of the exact fractional differintegrator and the approximation one can be obtained as Fig. \ref{Fig 6}.

\begin{figure}[!htbp]
  \centering
  \includegraphics[width=1.0\textwidth]{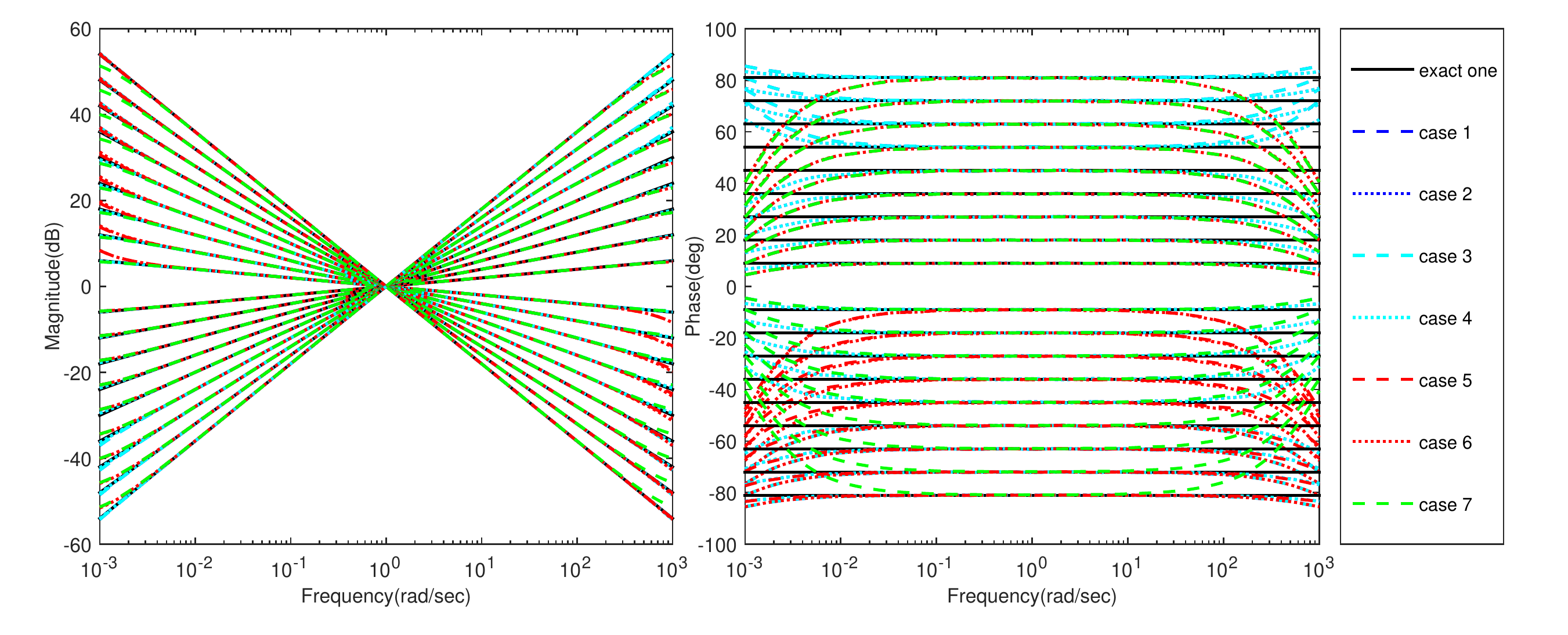}
  \caption{Bode diagram of fractional differintegrator and its approximate models with different $\alpha$. (Left: magnitude-frequency characteristics; Right: phase-frequency characteristics)}\label{Fig 6}
\end{figure}

Fig. \ref{Fig 6} clearly demonstrates that all the approximate models $\{{\hat {\mathscr I}_k^\alpha }\left( s \right),{\hat {\mathscr D}_k^\alpha }\left( s \right)\}$ perform well on the magnitude characteristics. Note that a comparable performance on the phase characteristic can also be reached in the middle frequency range, while the error increases slightly at the edge of the frequency range. For evaluating the approximation further, two quantitative indexes are defined
\begin{eqnarray}\label{Eq27}
{\textstyle
 \left\{\begin{array}{rl}
 {E_M}(k) =&\hspace{-0pt} 20\lg | {{{\mathscr I}^\alpha }\left( {{\rm{j}}\omega(k) } \right)} |-20\lg | {{\hat{ \mathscr I}_k^\alpha }\left( {{\rm{j}}\omega(k) } \right)} |,\\
 {E_P}(k) =&\hspace{-0pt} \frac{{180}}{\pi }\arg {{\mathscr I}^\alpha }\left( {{\rm{j}}\omega(k) } \right)-\frac{{180}}{\pi }\arg {\hat{\mathscr I}_k^\alpha }\left( {{\rm{j}}\omega(k) } \right). \end{array}\right.
}\hspace{-6pt}
\end{eqnarray}
If the approximation of ${{\mathscr D}^\alpha }(s)$ is considered, the corresponding ${E_M}(k)$ and ${E_P}(k)$ can be defined like (\ref{Eq27}).

Table \ref{Table 2} and Table \ref{Table 3} show the corresponding indexes implied in Fig. \ref{Fig 6}, from which four observations can be found. (a) For the same case, the results in the two tables are the same except case 5, since $\{{\hat {\mathscr I}_5^\alpha }\left( s \right)$ and ${\hat {\mathscr D}_5^\alpha }\left( s \right)\}$ were proposed by different scholars. In theory, the data of case 6 are equal in Table \ref{Table 2} and Table \ref{Table 3} by chance. (b) cases 1 - 4 perform better than cases 5 - 7, which clearly demonstrate the superiority of the developed methods. (c) With careful selection of $\varepsilon$, Algorithms \ref{Algorithm 3} and \ref{Algorithm 4} reduce to Algorithms \ref{Algorithm 1} and \ref{Algorithm 2}, respectively. Therefore, they bring the same approximation error. (d) Compared with case 1, case 2 is better, since it has large range of pole-zero.

\begin{table}[!htbp]
\caption{The maximum approximation error for ${\mathscr I}^\alpha(s)$ with $\alpha=0.1,0.2,\cdots,0.9$.}\label{Table 2}
\centering
\scriptsize
\begin{tabular}{lrrrr}
\toprule
{}&$\|{E_M}\|_\infty$&$\|{E_M}\|_2$&$\|{E_P}\|_\infty$&$\|{E_P}\|_2$\\
\hline
case 1&1.3179&26.7095&0.0226$\times10^{3}$&0.6570$\times10^{3}$\\
case 2&\textbf{0.4533}&\textbf{9.7653}&\textbf{0.0140$\times10^{3}$}&\textbf{0.3909$\times10^{3}$}\\
case 3&1.3179&26.7095&0.0226$\times10^{3}$&0.6570$\times10^{3}$\\
case 4&\textbf{0.4533}&\textbf{9.7653}&\textbf{0.0140$\times10^{3}$}&\textbf{0.3909$\times10^{3}$}\\
case 5&2.3792&49.7195&0.0387$\times10^{3}$&1.1463$\times10^{3}$\\
case 6&2.4251&49.5543&0.0406$\times10^{3}$&1.1886$\times10^{3}$\\
case 7&2.6441&55.8089&0.0405$\times10^{3}$&1.2137$\times10^{3}$\\
\bottomrule
\end{tabular}
\end{table}
\begin{table}[!htbp]
\caption{The maximum approximation error for ${\mathscr D}^\alpha(s)$ with $\alpha=0.1,0.2,\cdots,0.9$.}\label{Table 3}
\centering
\scriptsize
\begin{tabular}{lrrrr}
\toprule
{}&$\|{E_M}\|_\infty$&$\|{E_M}\|_2$&$\|{E_P}\|_\infty$&$\|{E_P}\|_2$\\
\hline
case 1&1.3179&26.7095&0.0226$\times10^{3}$&0.6570$\times10^{3}$\\
case 2&\textbf{0.4533}&\textbf{9.7653}&\textbf{0.0140$\times10^{3}$}&\textbf{0.3909$\times10^{3}$}\\
case 3&1.3179&26.7095&0.0226$\times10^{3}$&0.6570$\times10^{3}$\\
case 4&\textbf{0.4533}&\textbf{9.7653}&\textbf{0.0140$\times10^{3}$}&\textbf{0.3909$\times10^{3}$}\\
case 5&2.6441&55.8089&0.0405$\times10^{3}$&1.2137$\times10^{3}$\\
case 6&2.4251&49.5543&0.0406$\times10^{3}$&1.1886$\times10^{3}$\\
case 7&2.6441&55.8089&0.0405$\times10^{3}$&1.2137$\times10^{3}$\\
\bottomrule
\end{tabular}
\end{table}
\end{example}

\begin{example}Performance in time domain.

Setting $u(t)=\sin(t)$, then one has
\begin{equation}\label{Eq28}
{\textstyle
\left\{\begin{array}{rl}
 x=&\hspace{-0pt}{\mathscr I}^\alpha{\mathscr I}^{1-\alpha}[u](t)=1-\cos(t),\\
 y=&\hspace{-0pt}{\mathscr D}^{\alpha}{\mathscr I}^\alpha[u](t)=\sin(t),\\
 z=&\hspace{-0pt}{\mathscr D}^\alpha{\mathscr D}^{1-\alpha}[u](t)=\cos(t).
 \end{array}\right.}
\end{equation}
Changing $\{{ {\mathscr I}^\alpha }\left( s \right),{{\mathscr D}^\alpha }\left( s \right)\}$ to $\{{\hat {\mathscr I}_k^\alpha }\left( s \right),{\hat {\mathscr D}_k^\alpha }\left( s \right)\}$ yields the approximation signals $\hat x$, $\hat y$ and $\hat z$. The curves of these signals with $\alpha=0.4$ are plotted in Fig. \ref{Fig 7}.

\begin{figure}[!htbp]
  \centering
  \includegraphics[width=0.85\textwidth]{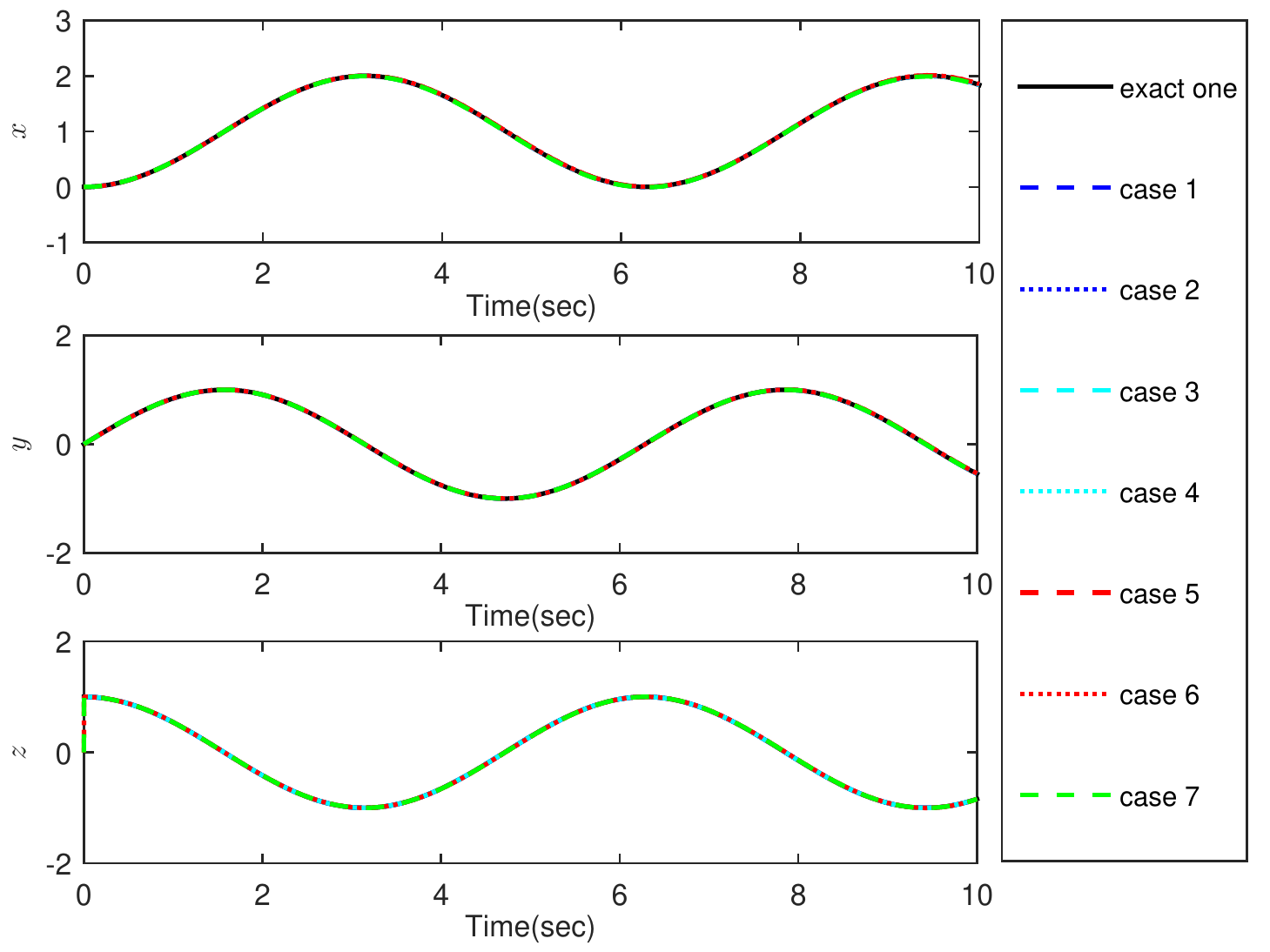}
  \caption{The approximation performance for three operations. (Above: $x$; Middle: $y$; Below: $z$)}\label{Fig 7}
\end{figure}

It is observed that these curves coincide with the exact one well. If we define the approximation error as $E_x=x-\hat x$, $E_y=y-\hat y$ and $E_z=z-\hat z$, then the results in Table \ref{Table 4} can be obtained. It is worth emphasizing that zero error approximation is achieved for the considered problem in cases 1 - 4. This is consistent with the results in Table \ref{Table 1}. The results in cases 5 - 7 also match Table \ref{Table 1} well. 

\begin{table}[!htbp]
\caption{The approximation error in time domain with $\alpha=0.4$.}\label{Table 4}
\centering
\scriptsize
\begin{tabular}{lcccccc}
\toprule
{}&$\|{E_x}\|_\infty$&$\|{E_x}\|_2$&$\|{E_y}\|_\infty$&$\|{E_y}\|_2$&$\|{E_z}\|_\infty$&$\|{E_z}\|_2$\\
\hline
case 1&0&0&0&0&0&0\\
case 2&0&0&0&0&0&0\\
case 3&0&0&0&0&0&0\\
case 4&0&0&0&0&0&0\\
case 5&0.0113&0.3880&0.0064&0.3063&1.0000&1.1004\\
case 6&0.0145&0.5289&0.0077&0.3695&1.0000&1.1150\\
case 7&0.0133&0.5277&0&0&1.0000&1.1004\\
\bottomrule
\end{tabular}
\end{table}

$\alpha=0.5$ is the singular case of this piecewise model and the related results are shown in Table \ref{Table 5}. Both ${\hat{ \mathscr I}_\kappa^\alpha }\left( s \right){\hat{ \mathscr I}_\kappa^{1 - \alpha }}\left( s \right) =\frac{1}{s}$ and ${\hat{ \mathscr D}_\kappa^\alpha }\left( s \right){\hat { \mathscr D}_\kappa^{1-\alpha }}\left( s \right) = s$ are not available in cases 1 - 4 which is unexpected. Thus, how to solve this problem will be an interesting subject in future.
\begin{table}[!htbp]
\caption{The approximation error in time domain with $\alpha=0.5$.}\label{Table 5}
\centering
\scriptsize
\begin{tabular}{lcccccc}
\toprule
{}&$\|{E_x}\|_\infty$&$\|{E_x}\|_2$&$\|{E_y}\|_\infty$&$\|{E_y}\|_2$&$\|{E_z}\|_\infty$&$\|{E_z}\|_2$\\
\hline
case 1&0.0146&0.5285&0&0&1.0000&1.1157\\
case 2&0.0126&0.4199&0&0&1.0000&1.1055\\
case 3&0.0146&0.5285&0&0&1.0000&1.1157\\
case 4&0.0126&0.4199&0&0&1.0000&1.1055\\
case 5&0.0114&0.3849&0.0067&0.3233&1.0000&1.1018\\
case 6&0.0146&0.5291&0.0078&0.3721&1.0000&1.1157\\
case 7&0.0135&0.5271&0&0&1.0000&1.1018\\
\bottomrule
\end{tabular}
\end{table}
\end{example}

\section{Conclusions}\label{Section 4}
In this paper, we firstly discussed the conflict on the order of the approximate model, i.e., ${\hat{ \mathscr I}_\kappa^\alpha }\left( s \right){\hat{ \mathscr I}_\kappa^{1 - \alpha }}\left( s \right) \neq \frac{1}{s}$, ${\hat{ \mathscr D}_\kappa^\alpha }\left( s \right){\hat { \mathscr I}_\kappa^{\alpha }}\left( s \right) \neq 1$ and ${\hat{ \mathscr D}_\kappa^\alpha }\left( s \right){\hat{ \mathscr D}_\kappa^{1 - \alpha }}\left( s \right) \neq s$ for $\alpha\in(0,1)$. Bearing this in mind, a piecewise model was proposed, which was able to avoid the conflict. Afterwards, two kinds of parameter design schemes were provided. From this, a detailed discussion was given on the numerical simulation and the circuit realization, etc. Finally, by means of comparison from the frequency domain and the time domain, the superiority of our methods in conflict avoidance and approximation accuracy were clearly illustrated.

\begin{acknowledgment}
The work described in this paper was fully supported by the National Natural Science Foundation of China (Grant No. 61601431; Funder ID: 10.13039/501100001809), the Anhui Provincial Natural Science Foundation (Grant No. 1708085QF141; Funder ID: 10.13039/501100003995), the Fundamental Research Funds for the Central Universities (Grant No. WK2100100028; Funder ID: 10.13039/501100012226), the General Financial Grant from the China Postdoctoral Science Foundation (Grant No. 2016M602032; Funder ID: 10.13039/501100002858) and the fund of China Scholarship Council (Grant No. 201806345002; Funder ID: 10.13039/501100004543).
\end{acknowledgment}

%

%
% ...

\bibliographystyle{asmems4}

\bibliography{database}

\begin{thebibliography}{10}

\bibitem{Sun:2018CNSNS}
Sun, H.~G., Zhang, Y., Baleanu, D., Chen, W., and Chen, Y.~Q., 2018.
\newblock ``A new collection of real world applications of fractional calculus
  in science and engineering''.
\newblock {\em Communications in Nonlinear Science and Numerical Simulation,
  {\bf 64}}, pp.~213--231.

\bibitem{Chen:2019NN}
Chen, L.~P., Huang, T.~w., Machado, J. A.~T., Lopes, A.~M., and Wu, R.~C.,
  2019.
\newblock ``Delay-dependent criterion for asymptotic stability of a class of
  fractional-order memristive neural networks with time-varying delays''.
\newblock {\em Neural Networks, {\bf 118}}, pp.~289--299.

\bibitem{Semary:2019JAR}
Semary, M.~S., Fouda, M.~E., Hassan, H.~N., and Radwan, A.~G., 2019.
\newblock ``Realization of fractional-order capacitor based on passive
  symmetric network''.
\newblock {\em Journal of Advanced Research, {\bf 18}}, pp.~147--159.

\bibitem{Adhikary:2018TCSI}
Adhikary, A., Choudhary, S., and Sen, S., 2018.
\newblock ``Optimal design for realizing a grounded fractional order inductor
  using {GIC}''.
\newblock {\em IEEE Transactions on Circuits and Systems I: Regular Papers,
  {\bf 65}}(8), pp.~2411--2421.

\bibitem{Adhikary:2016TCSI}
Adhikary, A., Sen, S., and Biswas, K., 2016.
\newblock ``Practical realization of tunable fractional order parallel
  resonator and fractional order filters''.
\newblock {\em IEEE Transactions on Circuits and Systems I: Regular Papers,
  {\bf 63}}(8), pp.~1142--1151.

\bibitem{Alejandro:2020JAR}
Silva-Ju\'{a}rez, A., Tlelo-Cuautle, E., Fraga, L. G. D.~L., and Li, R., 2020.
\newblock ``{FPGA}-based implementation of fractional-order chaotic oscillators
  using first-order active filter blocks''.
\newblock {\em Journal of Advanced Research, {\bf 25}}, pp.~77--85.

\bibitem{Zhang:2017ISA}
Zhang, X.~F., and Chen, Y.~Q., 2017.
\newblock ``Admissibility and robust stabilization of continuous linear
  singular fractional order systems with the fractional order $\alpha$: the
  $0<\alpha<1$ case''.
\newblock {\em ISA Transactions, {\bf 82}}, pp.~42--50.

\bibitem{Zhang:2020ISA}
Zhang, X.~F., and Wang, Z., 2020.
\newblock ``Stability and robust stabilization of uncertain switched fractional
  order systems''.
\newblock {\em ISA Transactions, {\bf 103}}, pp.~1--9.

\bibitem{Chen:2020TCSII}
Chen, L.~P., Wu, R.~C., Cheng, Y., and Chen, Y.~Q., 2020.
\newblock ``Delay-dependent and order-dependent stability and stabilization of
  fractional-order linear systems with time-varying delay''.
\newblock {\em IEEE Transactions on Circuits and Systems II: Express Briefs,
  {\bf 67}}(6), pp.~1064--1068.

\bibitem{Shi:2020VRa}
Shi, R.~Q., Li, Y., and Wang, C.~H., 2020.
\newblock ``Stability analysis and optimal control of a fractional-ordermodel
  for {African} swine fever''.
\newblock {\em Virus Research, {\bf 288}}.
\newblock doi: 198111.

\bibitem{Shi:2020VRb}
Shi, R.~Q., Lu, T., and Wang, C.~H., 2020.
\newblock ``Dynamic analysis of a fractional-order delayed model for hepatitis
  b virus with ctl immune response''.
\newblock {\em Virus Research, {\bf 277}}.
\newblock doi: 197841.

\bibitem{Vinagre:2000FCAA}
Vinagre, B.~M., Podlubny, I., Hernandez, A., and Feliu, V., 2000.
\newblock ``Some approximations of fractional order operators used in control
  theory and applications''.
\newblock {\em Fractional Calculus and Applied Analysis, {\bf 3}}(3),
  pp.~231--248.

\bibitem{Oustaloup:2000TCSI}
Oustaloup, A., Levron, F., Mathieu, B., and Nanot, F.~M., 2000.
\newblock ``Frequency-band complex noninteger differentiator: characterization
  and synthesis''.
\newblock {\em IEEE Transactions on Circuits and Systems I: Fundamental Theory
  and Applications, {\bf 47}}(1), pp.~25--39.

\bibitem{Poinot:2003SP}
Poinot, T., and Trigeassou, J.~C., 2003.
\newblock ``A method for modelling and simulation of fractional systems''.
\newblock {\em Signal Processing, {\bf 83}}(11), pp.~2319--2333.

\bibitem{Krishna:2011SP}
Krishna, B.~T., 2011.
\newblock ``Studies on fractional order differentiators and integrators: a
  survey''.
\newblock {\em Signal Processing, {\bf 91}}(3), pp.~386--426.

\bibitem{Meng:2012DSMC}
Meng, L., and Xue, D.~Y., 2012.
\newblock ``A new approximation algorithm of fractional order system models
  based optimization''.
\newblock {\em Journal of Dynamic Systems Measurement and Control, {\bf
  134}}(4).
\newblock id: 044504.

\bibitem{Romero:2013ISA}
Romero, M., De~Madrid, A.~P., Ma\~{n}oso, C., and Vinagre, B.~M., 2013.
\newblock ``{IIR} approximations to the fractional differentiator/integrator
  using {Chebyshev} polynomials theory''.
\newblock {\em ISA Transactions, {\bf 52}}(4), pp.~461--468.

\bibitem{Wei:2014IJCAS}
Wei, Y.~H., Gao, Q., Peng, C., and Wang, Y., 2014.
\newblock ``A rational approximate method to fractional order systems''.
\newblock {\em International Journal of Control, Automation and Systems, {\bf
  12}}(6), pp.~1180--1186.

\bibitem{Pakhira:2015ISA}
Pakhira, A., Das, S., Pan, I., and Das, S., 2015.
\newblock ``Symbolic representation for analog realization of a family of
  fractional order controller structures via continued fraction expansion''.
\newblock {\em ISA Transactions, {\bf 57}}, pp.~390--402.

\bibitem{Wei:2016ISA}
Wei, Y.~H., Tse, P.~W., Du, B., and Wang, Y., 2016.
\newblock ``An innovative fixed-pole numerical approximation for fractional
  order systems''.
\newblock {\em ISA Transactions, {\bf 62}}, pp.~94--102.

\bibitem{Tavazoei:2016JAS}
Tavazoei, M.~S., 2016.
\newblock ``Criteria for response monotonicity preserving in approximation of
  fractional order systems''.
\newblock {\em IEEE/CAA Journal of Automatica Sinica, {\bf 3}}(4),
  pp.~422--429.

\bibitem{Abdelaty:2018TCSII}
Abdelaty, A.~M., Elwakil, A.~S., Radwan, A.~G., Psychalinos, C., and Maundy,
  B.~J., 2018.
\newblock ``Approximation of the fractional-order {Laplacian} $s^\alpha$ as a
  weighted sum of first-order high-pass filters''.
\newblock {\em IEEE Transactions on Circuits and Systems II: Express Briefs,
  {\bf 65}}(8), pp.~1114--1118.

\bibitem{Deniz:2016ISA}
Deniz, F.~N., Alagoz, B.~B., Tan, N., and Atherton, D.~P., 2016.
\newblock ``An integer order approximation method based on stability boundary
  locus for fractional order derivative/integrator operators''.
\newblock {\em ISA Transactions, {\bf 62}}, pp.~154--163.

\bibitem{Sabatier:2018Alg}
Sabatier, J., 2018.
\newblock ``Solutions to the sub-optimality and stability issues of recursive
  pole and zero distribution algorithms for the approximation of fractional
  order models''.
\newblock {\em Algorithms, {\bf 11}}(7).
\newblock id. 103.

\bibitem{De:2018ISA}
De~Keyser, R., Muresan, C.~I., and Ionescu, C.~M., 2018.
\newblock ``An efficient algorithm for low-order direct discrete-time
  implementation of fractional order transfer functions''.
\newblock {\em ISA Transactions, {\bf 74}}, pp.~229--238.

\bibitem{Bai:2018ISA}
Bai, L., and Xue, D.~Y., 2018.
\newblock ``Universal block diagram based modeling and simulation schemes for
  fractional-order control systems''.
\newblock {\em ISA Transactions, {\bf 82}}, pp.~153--162.

\bibitem{Liang:2014IJSS}
Liang, S., Peng, C., Liao, Z., and Wang, Y., 2014.
\newblock ``State space approximation for general fractional order dynamic
  systems''.
\newblock {\em International Journal of Systems Science, {\bf 45}}(10),
  pp.~2203--2212.

\bibitem{Liang:2017JCSC}
Liang, Y.~S., and Lu, J.~L., 2017.
\newblock ``Direct low order rational approximations for fractional order
  systems in narrow frequency band: a fix-pole method''.
\newblock {\em Journal of Circuits, Systems and Computers, {\bf 26}}(04).
\newblock id. 1750065.

\bibitem{Wei:2018ISA}
Wei, Y.~H., Wang, J.~C., Liu, T.~Y., and Wang, Y., 2019.
\newblock ``Fixed pole based modeling and simulation schemes for fractional
  order systems''.
\newblock {\em ISA Transactions, {\bf 84}}, pp.~43--54.

\bibitem{Li:2020ISA}
Li, A., Wei, Y.~H., and Wang, Y., 2020.
\newblock ``A numerical approximation method for fractional order systems with
  new distributions of zeros and poles''.
\newblock {\em ISA Transactions, {\bf 99}}, pp.~20--27.

\bibitem{Wei:2021DSMC}
Wei, Y.~H., Zhang, H., Hou, Y.~Q., and Cheng, K., 2021.
\newblock ``Multiple fixed pole based rational approximation for fractional
  order systems''.
\newblock {\em Journal of Dynamic Systems Measurement and Control}.
\newblock doi: 10.1115/1.4049557.

\end{thebibliography}

%

\end{document}